\documentclass[11pt,reqno]{article}
\usepackage{euscript,amstext,amsthm,amssymb}
\usepackage{amsmath}

\newcommand{\be}{\begin{equation}}
      \newcommand{\ee}{\end{equation}}
      \newcommand{\ba}{\begin{eqnarray}}
       \newcommand{\ea}{\end{eqnarray}}
\newcommand{\ban}{\begin{eqnarray*}}
       \newcommand{\ean}{\end{eqnarray*}}

\newcommand{\pt}{\partial}

\newcommand{\ra}{\rightarrow}

 \renewcommand{\o}[2]{\frac{#1}{#2}}

\newcommand{\bM}{\partial M}

\newcommand{\Mas}{\mbox{\rm Mas}}

\newcommand{\ind}{\mbox{\rm ind}}

\begin{document}
\newtheorem{lem}[theo]{Lemma}
\newtheorem{prop}[theo]{Proposition}
\newtheorem{coro}[theo]{Corollary}

\title{An Index Theorem for Toeplitz Operators\\ on Odd Dimensional Manifolds with Boundary}
\author{Xianzhe Dai\thanks{Partially supported by NSF.}$\ $ and Weiping
Zhang\thanks{Partially supported by MOEC, MOSTC and NNSFC.}}
%\date{\it Dedicated to our teachers Jean-Michel Bismut and Jeff Cheeger}

\maketitle
\begin{abstract}
We establish an index theorem for Toeplitz operators on odd
dimensional spin manifolds with boundary. It may be thought of as
an odd dimensional analogue of the Atiyah-Patodi-Singer index
theorem for Dirac operators on manifolds with boundary. In
particular, there occurs naturally an invariant of $\eta$ type
associated to  $K^1$ representatives on  even dimensional
manifolds, which should be of independent interests. For example,
it gives an intrinsic interpretation of the so called Wess-Zumino
term in the WZW theory in physics.
\end{abstract}

%\noindent
\section{Introduction}

%$\ $

On an even dimensional  smooth closed spin Riemannian manifold
$M$, let $S(TM)$ be the corresponding bundle of spinors over $M$
and $E$ be a Hermitian vector bundle over $M$ equipped with a
Hermitian connection.  The (twisted) Dirac operator
$D^E:\Gamma(S(TM)\otimes E) \rightarrow \Gamma(S(TM)\otimes E) $
is elliptic and self-adjoint. Since $\dim M$ is even, the spinors
split:
$$ S(TM)\otimes E = S^{+}(TM)\otimes E \oplus S^{-}(TM)\otimes E, $$
in terms of which the Dirac operator is off diagonal:
$$ D^E = \left( \begin{array}{cc} 0 & D^E_- \\
D^E_+ & 0 \end{array} \right). $$ The Atiyah-Singer index theorem
expresses the index of $D^E_+$ in terms of the characteristic
numbers:
$$ \ind \, D^E_+ = \left\langle \widehat{A}(TM){\rm ch}(E),[M]\right\rangle $$
where $\widehat{A}(TM)$ is the Hirzebruch $\widehat{A}$-class of
$TM$, ${\rm ch}(E)$ is the Chern character of $E$ (cf. [Z1, Chap.
1]).

Now let $M$ be an {\em odd} dimensional  smooth closed spin
Riemannian manifold. Any elliptic differential operator on $M$ will
have index zero. In this case, the appropriate index to consider is
that of Toeplitz operators. It also fits perfectly with the
interpretation of the index of Dirac operator on even dimensional
manifolds as a pairing between the even $K$-group and $K$-homology.
Thus in the odd dimensional case one considers the odd $K$-group and
odd $K$-homology. An element of $K^{-1}(M)$ can be  represented by a
differentiable map from $M$ into the general linear group
$$ g: \ M \longrightarrow { GL}(N, {\bf C}), $$
where $N$ is a positive integer.
As we mentioned the appropriate index pairing between the odd $K$-group
and $K$-homology is given by that of the Toeplitz
operator, defined as follows.

First of all, $L^2(S(TM)\otimes E)$, the natural $L^2$-completion of $\Gamma(S(TM)\otimes E) $,
 splits into an orthogonal
direct sum as
$$L^2(S(TM)\otimes E)=\bigoplus_{\lambda\in {\rm Spec}(D^E)} E_\lambda,$$
where $E_\lambda$ is the eigenspace associated to the eigenvalue  $\lambda$ of $D^E$. Set
$$L^2_+(S(TM)\otimes E) =\bigoplus_{\lambda\geq 0}E_\lambda, $$
and denote by $P^E_{\geq 0}$
the orthogonal projection from $L^2(S(TM)\otimes E)$ to $L^2_+(S(TM)\otimes E)$.

Now consider the trivial vector bundle ${\bf C}^N$ over $M$.
We equip ${\bf C}^N$ with the canonical trivial metric and connection. Then $P^E_{\geq 0}$ extends naturally
to an orthogonal projection from $L^2(S(TM)\otimes E \otimes {\bf C}^N)$ to
$L^2_+(S(TM)\otimes E \otimes {\bf C}^N)$ by acting
as identity on ${\bf C}^N$. We still denote
this extension by $P^E_{\geq 0}$.

The map $g$ can be interpreted as automorphism of the trivial complex vector bundle ${\bf C}^N$.
%$g\in {\rm Aut}({\bf C}^N)$.
 Moreover $g$
extends naturally to an action on $L^2(S(TM)\otimes E\otimes {\bf C}^N)$ by acting as identity on
$L^2(S(TM)\otimes E)$. We still denote this extended action by $g$.

With the above data given, the Toeplitz operator $T^E_g$ can be defined as
$$T^E_g=P^E_{\geq 0}gP^E_{\geq 0}:L^2_+\left(S(TM)\otimes E \otimes {\bf C}^N\right)
\longrightarrow L^2_+\left(S(TM)\otimes E \otimes {\bf C}^N\right).$$

The first important fact is that $T^E_g$ is a Fredholm operator. Moreover, it is equivalent to an elliptic
pseudodifferential operator of order zero. Thus one can compute its index by using the
Atiyah-Singer index theorem [AS], as was indicated in the paper of Baum and Douglas [BD]:
$${\rm ind}\, T^E_g=-\left\langle \widehat{A}(TM){\rm ch}(E){\rm ch}(g),[M]\right\rangle ,\eqno(1.1)$$
where ${\rm ch}(g)$ is the odd Chern character associated to $g$
(cf. [Z1, Chap. 1]).

There is also an analytic proof of (1.1) by using  heat kernels. For this one first
note that by a simple deformation, one may well assume that $g$ is unitary. Then a result of
Booss and Wojciechowski (cf. [BW]) shows that the computation of ${\rm ind}\, T^E_g$ is equivalent
to the computation of the spectral flow of the linear family of self-adjoint elliptic operators,
acting of $\Gamma(S(TM)\otimes E\otimes {\bf C}^N)$, which connects $D^E$ and $gD^Eg^{-1} $.
The resulting spectral flow can then be computed by  variations of $\eta$-invariants, where the
heat kernels are naturally involved.
These ideas have been extended in [DZ] to give a heat kernel proof of a family extension of (1.1).

The purpose of this paper is to establish a generalization of
(1.1) to the case where $M$ is a spin manifold with boundary
$\partial M$, by extending the above heat kernel proof strategy.
We wish to point out that when $g|_{\partial M}$ is the identity,
such a generalization can be reduced easily to a result of
Douglas and Wojciechowski [DW]. Thus the main concern for us in
this paper will be the case where $g|_{\partial M}$ is not the
identity.

A full statement of our main result will be given in Section 2
(Theorem 2.3). Here we only point out that our formula may be
viewed as an odd dimensional analogue of the Atiyah-Patodi-Singer
index theorem [APS1] for Dirac operators on even dimensional
manifolds with boundary. In particular, a very interesting
invariant of $\eta$-type for even dimensional manifolds and $K^1$
representaitves appears in our formula, which plays a role similar
to that played by the $\eta$ invariant term in the
Atiyah-Patodi-Singer index theorem. There is also an interesting
new integer term here, a triple Maslov index introduced in [KL].

This paper is organized as follows. In Section 2, we introduce the
notations and state the main result of this paper. In Section 3,
we introduce a perturbation to overcome a technical difficulty and
prove an index formula for the perturbed Toeplitz operator. In
Section 4, we compare the index of Toeplitz operator and that of
the perturbed one and prove our main result. In Section 5 we
discuss some generalizations of the main result proved in Sections
3, 4, including the basic properties of the $\eta$-type invariant
mentioned above. We also include an Appendix in which we outline a
new proof of (1.1) by a simple use of the Atiyah-Patodi-Singer
index theorem.

$\ $

\noindent {\bf Acknowledgments}. The authors wish to thank
Professor Jerry Kaminker who suggested to them the problem
considered in this paper. They also acknowledge interesting
discussions with Professors Jerry Kaminker and Henri Moscovici.
Part of the work of the second author was done while he was
visiting MSRI during the program of Spectral Invariants, and
during his visit to MIT for the Spring semester of 2001. He would
like to thank the organizers (in particular, Professor Alice S.-Y.
Chang) of the Spectral Invariants program for invitation, and MSRI
for financial support. He is also grateful to Professors Richard
Melrose and Gang Tian for arranging his visit to MIT, and to MIT
for financial support. Part of the revision was done while the
first author was visiting the Nankai Institute (now the Chern
Institute) of Mathematics in July, 2005. He would like to thank
the Chern Institute of Mathematics for hospitality.  Finally, the
authors thank the referee who found a gap in an earlier version
and made many useful suggestions.

%\noindent { \bf \S 1.
\section{An index theorem for Toeplitz operators on manifolds with boundary}

%$\ $

In this section, we state the main result of this paper, which extends (1.1) to manifolds
with boundary.

This section is organized as follows. In a), we present our basic
geometric data and define the Toeplitz operators on manifolds with
boundary. In b), we define an $\eta$-type invariant for $K^1$
representatives on even dimensional manifolds, which will appear
in the statement of the main result. In c), we state the main
result of this paper, the proof of which will be presented in the
next two sections.

$\ $

\noindent {\bf a). Toeplitz operators on manifolds with boundary}

$\ $

Let $M$ be an odd dimensional oriented spin manifold with boundary
$\partial M$. We assume that $M$ carries a fixed spin structure.
Then $\partial M$ carries the canonically induced orientation and
spin structure. Let $g^{TM}$ be a Riemannian metric on $TM$ such
that it is of product structure near the boundary $\partial M$. That
is, there is a tubular neighborhood, which, without loss of
generality, can be taken to be $[0, 1)\times
\partial M \subset M$ with $\partial M=\{
0\}\times
\partial M$ such that
$$\left.g^{TM}\right|_{[0, 1)\times \partial M}=dx^2\oplus g^{T\partial M},\eqno(2.1)$$
where $x\in [0, 1)$ is the geodesic distance to $\bM$ and
$g^{T\partial M}$ is the restriction of $g^{TM}$ on $\partial M$.
Let $\nabla^{TM}$ be the Levi-Civita connection of $g^{TM}$. Let
$S(TM)$ be the Hermitian bundle of spinors associated to
$(M,g^{TM})$. Then $\nabla^{TM}$ extends naturally to a Hermitian
connection $\nabla^{S(TM)}$ on $S(TM)$.

Let $E$ be a Hermitian vector bundle over $M$. Let $\nabla^E$ be a
Hermitian connection on $E$. We assume that the Hermitian metric
$g^E$ on $E$ and connection $\nabla^E$ are of product structure over
$[0, 1)\times \partial M$. That is, if we denote $\pi:[0, 1)\times
\partial M \rightarrow \partial M$ the natural projection, then
$$
\left. g^E\right|_{[0, 1)\times \partial M}=\pi^* \left(\left. g^E\right|_{\partial M}\right),
\ \ \left. \nabla^E\right|_{[0, 1)\times \partial M}= \pi^*
\left(\left. \nabla^E\right|_{\partial M}\right).\eqno(2.2)
$$

For any $X\in TM$, we extend the Clifford action $c(X)$ of $X$ on $S(TM)$ to an action on
$S(TM)\otimes E$ by acting as identity on $E$, and still denote this extended action by $c(X)$.
Let $\nabla^{S(TM)\otimes E}$ be the tensor product connection on $S(TM)\otimes E$ obtained
from $\nabla^{S(TM)}$ and $\nabla^E$.

The canonical
(twisted) Dirac operator $D^E$ is defined by
$$D^E=\sum_{i=1}^{\dim M}c(e_i)\nabla^{S(TM)\otimes E}_{e_i}:\Gamma(S(TM)\otimes E)
\longrightarrow \Gamma(S(TM)\otimes E),  \eqno(2.3)$$ where
$e_1,\dots, e_{\dim M}$ is an orthonormal basis of $TM$. By (2.1)
and (2.2), over $[0, 1)\times \partial M$, one has
$$D^E=c\left({\partial\over \partial x}\right)\left({\partial\over\partial x}
+\pi^*D^E_{\partial M}\right),\eqno(2.4)$$ where $D^E_{\partial
M}:\Gamma((S(TM)\otimes E)|_{\partial M})\rightarrow
\Gamma((S(TM)\otimes E)|_{\partial M})$ is the induced Dirac
operator on $\partial M$.

We now introduce the APS type boundary conditions for $D^E$. The
induced Dirac operator on the boundary, $D^E_{\partial M}$, is
elliptic and self-adjoint. Let $L^2_{+}((S(TM) \otimes E)|_{\partial
M})$ be the space of the direct sum of eigenspaces of positive
eigenvalues of $D^E_{\partial M}$. Let $P_{\partial M}$ denote the
orthogonal projection operator from $L^2((S(TM)\otimes E)|_{\partial
M})$ to $L^2_{+}((S(TM)\otimes E)|_{\partial M})$ (for simplicity we
suppress the dependence on $E$).

As is well known, the APS projection $P_{\partial M}$ is an elliptic
global boundary condition for $D^E$. However, to get self adjoint
boundary conditions, we need to modify it by a Lagrangian subspace
of $\ker D^E_{\partial M}$, namely, a subspace $L$ of $\ker
D^E_{\partial M}$ such that $c({\partial \over \partial x})L=L^\perp
\cap (\ker D^E_{\partial M})$. Since $\partial M$ bounds $M$, by the
cobordism invariance of the index, such Lagrangian subspaces always
exist.

The modified APS projection is obtained by adding the projection
onto the Lagrangian subspace. Let $P_{\partial M}(L)$ denote the
orthogonal projection operator from $L^2((S(TM) \otimes
E)|_{\partial M})$ to $L^2_{+}((S(TM)\otimes E)|_{\partial
M})\oplus L$:
$$P_{\partial M}(L) =P_{\partial M}+P_L,\eqno(2.5)$$
where $P_L$ denotes the orthogonal projection  from
$L^2((S(TM)\otimes E)|_{\partial M})$ to $L$.

The pair $(D^E, P^E_{\partial M}(L) )$ forms a self-adjoint
elliptic boundary problem, and $P_{\partial M}(L) $ is called an
Atiyah-Patodi-Singer boundary condition associated to $L$. We will
also denote the corresponding elliptic self-adjoint operator by
$D^E_{P_{\partial M}(L) }$.

  Let $L^{2, +}_{P_{\partial M}(L)}  (S(TM)\otimes E))$
be the space of the direct sum of eigenspaces of {\it
non-negative} eigenvalues of $D^E_{P_{\partial M}(L) }$. This can
be viewed as an analog of the Hardy space. We denote by
$P_{P_{\partial M}(L)}$ the orthogonal projection from
$L^2(S(TM)\otimes E)$ to $L^{2, +}_{P_{\partial
M}(L)}(S(TM)\otimes E))$.

Let $N>0$ be a positive integer, let ${\bf C}^N$ be the trivial
complex vector bundle over $M$ of rank $N$, which carries the
trivial Hermitian metric and the trivial Hermitian connection.
Then all the above construction can be developed in the same way
if one replaces $E$ by $E\otimes {\bf C}^N$. And all the operators
considered here extend to act on ${\bf C}^N$ by identity. If there
is no confusion we will also denote them by their original
notation.

Now let $g:M\rightarrow GL(N,{\bf C})$ be a smooth automorphism of
${\bf C}^N$. With simple deformation, we can assume that $g$ is
unitary. That is, $g:M\rightarrow U(N)$. Furthermore, we make the
assumption that $g$ is of product structure over $[0, 1)\times
\partial M$, that is,
$$g|_{[0, 1)\times \partial M}=\pi^*\left(g|_{\partial M}\right).\eqno(2.6)$$

Clearly, $g$ extends to an action on $S(TM)\otimes E\otimes {\bf C}^N$ by acting as identity
on $S(TM)\otimes E$. We still denote this extended action by $g$.

Since $g$ is unitary, one verifies easily that the operator
$gP_{\partial M}(L) g^{-1}$ is again an orthogonal projection on
$L^2((S(TM)\otimes E\otimes {\bf C}^N)|_{\partial M})$, and that
$gP_{\partial M}(L) g^{-1}- P_{\partial M}(L) $ is a
pseudodifferential operator of order less than zero. Moreover, the
pair $(D^{E},gP_{\partial M}(L) g^{-1})$ forms a self-adjoint
elliptic boundary problem. We denote its associated elliptic
self-adjoint operator by $D^{E}_{gP_{\partial M}(L) g^{-1}}$. Thus
$D^{E}_{gP_{\partial M}(L) g^{-1}}$ has the boundary condition which
is the conjugation by $g$ of the previous APS type condition.

The necessity of using the conjugated boundary condition here is
from the fact that, if $s\in L^2(S(TM)\otimes E\otimes {\bf C}^N)$
verifies $P_{\partial M}(L) (s|_{\partial M})=0$, then $gs$ verifies
$gP_{\partial M}(L) g^{-1}((gs)|_{\partial M}) =0$.

Thus, consider also the analog of Hardy space for the conjugated
boundary value problem, $L^{2, +}_{gP_{\partial
M}(L)g^{-1}}(S(TM)\otimes E\otimes {\bf C}^N)$ which is the space of
the direct sum of eigenspaces of {\em nonnegative} eigenvalues of
$D^{E}_{gP_{\partial M}(L) g^{-1}}$. Let $P_{gP_{\partial M}(L)
g^{-1}}$ denote the orthogonal projection   from $L^2(S(TM)\otimes
E\otimes {\bf C}^N)$ to $L^{2, +}_{gP_{\partial M}(L) g^{-1}}
(S(TM)\otimes E\otimes {\bf C}^N)$.

$\ $

\noindent {\bf Definition 2.1.} {The Toeplitz operator $T^E_g(L)$ is defined by
$$T^E_g(L)=P_{gP_{\partial M}(L) g^{-1}} \circ
g \circ P_{P_{\partial M}(L)}:$$
$$L^{2, +}_{P_{\partial M}(L)}
\left(S(TM)\otimes E\otimes {\bf C}^N\right) \rightarrow L^{2,
+}_{gP_{\partial M}(L) g^{-1}} \left(S(TM)\otimes E\otimes {\bf
C}^N\right). \eqno(2.7)$$}

One verifies that $T^E_g(L)$ is a Fredholm operator. The main
purpose of this paper is to establish an index formula for it in
terms of geometric data.

$\ $

\noindent {\bf b). Perturbation}

$\ $

The analysis of the conjugated elliptic boundary value problem
$D^{E}_{gP_{\partial M}(L) g^{-1}}$ turns out to be surprisingly
subtle and difficult. To circumvent this difficulty, we now
construct a perturbation of the original problem.

Let $\psi=\psi(x)$ be a cut off function which is identically $1$ in
the $\epsilon$-tubular neighborhood of $\partial M$ ($\epsilon
>0$ sufficiently small) and vanishes outside the $2\epsilon$-tubular
neighborhood of $\partial M$. Consider the Dirac type operator
$$
D^{\psi}=(1-\psi)D^E + \psi gD^E g^{-1}. \eqno(2.8)
$$

The effect of this perturbation is that, near the boundary, the
operator $D^{\psi}$ is actually given by the conjugation of $D^E$,
and therefore, the elliptic boundary problem $(D^{\psi},
gP_{\partial M}(L) g^{-1})$ is now the conjugation of the APS
boundary problem $(D^E, P_{\partial M}(L))$.

All previous consideration applies to $(D^{\psi}, gP_{\partial M}(L)
g^{-1})$ and its associated self adjoint elliptic operator $D^{\psi}_{gP_{\partial M}(L)
g^{-1}}$. In particular, we have the perturbed Toeplitz operator
$$
T^E_{g, \psi}(L)=P^{\psi}_{gP_{\partial M}(L) g^{-1}} \circ
g \circ P_{P_{\partial M}(L)}:$$
$$ L^{2, +}_{P_{\partial M}(L)}
\left(S(TM)\otimes E\otimes {\bf C}^N\right) \rightarrow L^{2, +,
\psi}_{gP_{\partial M}(L) g^{-1}} \left(S(TM)\otimes E\otimes {\bf
C}^N\right), \eqno(2.9)
$$
where $P^{\psi}_{gP_{\partial M}(L) g^{-1}}$ is the APS projection
associated to $D^{\psi}_{gP_{\partial M}(L) g^{-1}}$, whose range is
denoted by $L^{2, +, \psi}_{gP_{\partial M}(L) g^{-1}}
\left(S(TM)\otimes E\otimes {\bf C}^N\right)$.

We will also need to consider the conjugation of $D^{\psi}$:
$$
D^{\psi, g}=g^{-1} D^{\psi} g =D^E + (1-\psi)
g^{-1}[D^E, g] .  \eqno(2.10)
$$

$\ $

\noindent {\bf c). An invariant of $\eta$ type for even dimensional
manifolds}

$\ $

Given an even dimensional closed spin manifold $X$, we consider
the cylinder $[0, 1]\times X$ with the product metric. Let $g: \ X
\ra U(N)$ be a map from $X$ into the unitary group which extends
trivially to the cylinder. Similarly, $E \ra X$ is an Hermitian
vector bundle which is also extended trivially to the cylinder. We
make the assumption that ${\rm ind}\, D^E_+=0$ on $X$.

Consider the analog of $D^{\psi,g}$ as defined in (2.10), but now on
the cylinder $[0, 1]\times X$  and denote it by $D^{\psi,g}_{[0,
1]}$. We equip it with the boundary condition $P_{X}(L) $ on one of
the boundary component $\{0\}\times X$ and the boundary condition
${\rm Id}-g^{-1}P_{X}(L) g $ on the other boundary component
$\{1\}\times X$ (Note that the Lagrangian subspace $L$ exists by our
assumption of vanishing index). Then $(D^{\psi,g}_{[0, 1]}, P_X(L) ,
{\rm Id}-g^{-1}P_{X}(L) g)$ forms a self-adjoint elliptic boundary
problem. For simplicity, we will still denote the corresponding
elliptic self-adjoint operator by $D^{\psi, g}_{[0, 1]}$.

Let $\eta(D^{\psi, g}_{[0, 1]},s)$ be the $\eta$-function of
$D^{\psi, g}_{[0, 1]}$ which, when ${\rm Re}(s)>>0$, is defined by
$$
\eta(D^{\psi, g}_{[0, 1]},s)
=\sum_{\lambda \neq 0}{{\rm sgn}(\lambda)\over
|\lambda|^s},\eqno(2.11)
$$
where $\lambda$ runs through the nonzero
eigenvalues of $D^{\psi, g}_{[0, 1]}$.

By \cite{dw,mu,df}, one knows that the $\eta$-function
$\eta(D^{\psi, g}_{[0, 1]},s)$ admits a meromorphic extension to ${\bf
C}$ with $s=0$ a regular point (and only simple poles). One then
defines, as in \cite{aps1}, the $\eta$-invariant of
$D^{\psi,g}_{[0,1]}$, denoted by
$\eta(D^{\psi,g}_{[0, 1]})$, to be the value at $s=0$ of
$\eta(D^{\psi,g}_{[0, 1]},s)$, and the reduced $\eta$-invariant by
{$$\overline{\eta}(D^{\psi,g}_{[0, 1]})={\dim \ker D^{\psi,g}_{[0, 1]} +
\eta(D^{\psi,g}_{[0, 1]}) \over 2}. \eqno(2.12)$$}

In our application, we will apply this construction to the cylinder
$[0, 1] \times \partial M$. i.e., $X=\partial M$ is a boundary. We
point out in passing that the invariant
$\overline{\eta}(D^{\psi,g}_{[0, a]})$, similarly constructed on a
cylinder $[0, a] \times X$,  does not depend on the radial size of
the cylinder $a>0$ by a rescaling argument (cf. [M\"u, Proposition
2.16]).
\newline

\noindent {\bf Definition 2.2.} {We define an invariant of $\eta$
type for the complex vector bundle $E$ on the even dimensional
manifold $X$ (with vanishing index) and the $K^1$ representative
$g$ by
$$
\overline{\eta}(X, g)=
 \overline{\eta}(D^{\psi, g}_{[0, 1]}) - {\rm sf} \left\{D^{\psi,g}_{[0, 1]}(s);
 0 \leq s \leq 1 \right\}, \eqno(2.13)
$$
where $D^{\psi,g}_{[0, 1]}(s)$ is a path connecting $g^{-1} D^E g$
with $D^{\psi,g}_{[0, 1]}$ defined by
$$
D^{\psi,g}(s)= D^E + (1-s\psi) g^{-1}[D^E, g]  \eqno(2.14)
$$
on $[0, 1]\times X$, with the boundary condition $P_{X}(L) $ on
$\{0\}\times X$ and the boundary condition ${\rm
Id}-g^{-1}P_{X}(L) g$ at $\{1\}\times X$.}
\newline

 We will show in Section 5 that $\overline{\eta}(X,
g)$ does not depend on the cut off function $\psi$.

$\ $

\noindent {\bf c). An index theorem for $T^E_g(L)$}

$\ $

Recall that $g: M \rightarrow U(N)$. Thus $g^{-1}dg$ defines a
$\Gamma({\rm End}({\bf C}^N))$-valued 1-form on $M$. Let ${\rm
ch}(g)$ denote the odd Chern character form of $g$ defined by (cf.
[Z1, Chap. 1])
$$
{\rm ch}(g)= \sum_{n=0}^{\dim M-1\over 2} {n!\over (2n+1)!}{\rm
Tr}\left[\left(g^{-1}dg\right)^{2n+1}\right]. \eqno(2.15)
$$

%For later purpose, we define similarly a perturbed version
%$$ {\rm ch}_{\psi}(g)=
%\int_0^1 {\rm Tr}\left[(1-\psi)g^{-1}dg\exp \left((1-u)(1-\psi) (u +
%\psi - u \psi) (g^{-1}dg)^2 \right) \right] du. $$

Recall also that $\nabla^{TM}$ is the Levi-Civita connection associated
to the Riemannian metric $g^{TM}$, and $\nabla^E$ is the Hermitian
connection on $E$. Let $R^{TM}=(\nabla^{TM})^2$ (resp.
$R^E=(\nabla^E)^2$) be the curvature of $\nabla^{TM}$ (resp.
$\nabla^E$).

Let ${\cal P}_M$ denote the Calder\'on projection associated to
$D^{E\otimes {\bf C}^N}$ on $M$ (cf. [BW]). Then ${\cal P}_M$ is an
orthogonal projection on $L^2((S(TM)\otimes E\otimes {\bf
C}^N)|_{\partial M})$, and that ${\cal P}_M -P_{\partial M}(L)$ is a
pseudodifferential operator of order less than zero.

Let $\tau_\mu (gP_{\partial M}(L) g^{-1},P_{\partial M}(L) ,
{\cal P}_M)\in {\bf Z}$ be the Maslov triple index in the sense of
Kirk and Lesch [KL, Definition 6.8]\footnote{Note the slight difference in
notation here: in KL's notation, our first variable should be $I-gP_{\partial M}(L) g^{-1}$}.

We can now state the main result of this paper as follows.

$\ $

\noindent {\bf Theorem 2.3.} {\it The following identity holds,
$$
 {\rm ind}\, T^E_g(L)=-\left({1\over 2\pi\sqrt{-1}}\right)^{(\dim M+1)/2}\int_M
\widehat{A}\left(R^{TM}\right){\rm Tr}\left[\exp\left(-R^E\right)\right]{\rm ch}(g,d) \
$$
$$
-  \overline{\eta}(\partial M, g) + \tau_\mu \left(gP_{\partial
M}(L) g^{-1},P_{\partial M}(L) ,  {\cal P}_M\right) .
\hspace{.2in} \eqno(2.16)
$$}

\noindent{\bf Remark 2.4}.  We will show in Theorem 4.3 that the
same formula holds without the product type assumption (2.6). Also,
the spin assumption can be relaxed and the same result holds for
general Dirac type operators, in particular, spin$^c$ Dirac
operators.

$\ $

\noindent{\bf Remark 2.5}. Our formula (2.16) is closely related
to the so called WZW theory in physics \cite{w}. When $\bM=S^2$ or
a compact Riemann surface and $E$ is trivial, the local term in
(2.16) is precisely the Wess-Zumino term, which allows an integer
ambiguity, in the WZW theory. Thus, our eta invariant
$\overline{\eta}(\partial M, g)$ gives an intrinsic interpretation
of the Wess-Zumino term without passing to the bounding
$3$-manifold. In fact, for $\bM=S^2$, it can be further reduced to
a local term on $S^2$ by using Bott's periodicity, see {\bf Remark
5.9}.

$\ $

The following immediate consequence is of independent interests and
will be studied further in Section 4.

$\ $

\noindent {\bf Corollary 2.6.} {\it The number
$$
\left({1\over 2\pi\sqrt{-1}}\right)^{(\dim M+1)/2}\int_M
\widehat{A}\left(R^{TM}\right){\rm
Tr}\left[\exp\left(-R^E\right)\right]{\rm ch}(g,d)+
\overline{\eta}(\partial M, g)
$$
is an integer.}

$\ $

The next two sections will be devoted to a proof of Theorem 2.3.

%$$\ $$

%\noindent {\bf \S 2.

\section{$\eta$-invariants, spectral flow and the index of \\ the perturbed Toeplitz operator}

%$\ $

In this section, we prove an index formula for the perturbed
Toeplitz operator $T^E_{g, \psi}(L)$. The strategy follows from that
of the heat kernel proof of (1.1) sketched in Introduction. However,
as we are dealing with the case of manifolds with boundary, we must
make necessary modifications  at each step of the procedure.

This section is organized as follows. In a), we reduce the
computation of ${\rm ind}\, T_{g, \psi}^E(L)$ to the computation of
a spectral flow of a natural family of self-adjoint elliptic
operators on manifolds with boundary. In b), we reduce the
computation of the above mentioned spectral flow to a computation of
certain $\eta$-invariants as well as their variations. We then apply
a result of Kirk and Lesch [KL, Theorem 7.7] to reduce the proof of
Theorem 2.3 to a computation of certain local index term arising
from the variations of $\eta$-invariants. In c), we prove the index
formula by computing the local index term.

$\ $

\noindent {\bf a). Spectral flow and the index of the perturbed
Toeplitz operators}

$\ $

Recall that $D^{\psi}$ defined in (2.8)  is the perturbed Dirac
operator on $M$ acting on $\Gamma(S(TM)\otimes E\otimes {\bf C}^N)$,
and $g:M\rightarrow U(N)$ is a smooth map.

For any $u\in [0,1]$, in view of (2.10), set
$$
D^{\psi,g}(u)=(1-u)D^E+ ug^{-1}D^{\psi}g=D^E +u(1-\psi)g^{-1}[D^E,
g].\eqno(3.1)
$$
Then  for each $u\in [0,1]$, the  boundary condition $P_{\partial
M}(L)$ is still a self-adjoint elliptic boundary condition for
$D^{\psi,g}(u)$. We denote the corresponding self-adjoint elliptic
operator by $D^{\psi, g}_{P_{\partial M}(L)}(u)$, which depends
smoothly on $u\in [0,1]$.

Let ${\rm sf}(D^{\psi, g}_{P_{\partial M}(L)}(u), 0\leq u \leq 1)$ be
the spectral flow of the this one parameter family of elliptic
self-adjoint operators in the sense of Atiyah-Patodi-Singer [APS2].

The following result generalizes a theorem of Booss-Wojciechowski
(cf. [BW, Theorem 17.17]) for closed manifolds.

$\ $

\noindent {\bf Theorem 3.1.} {\it We have,}
$${\rm ind}\, T^E_{g, \psi}(L)=-{\rm sf}\left(D^{\psi,g}_{P_{\partial M}(L)}(u), 0\leq u\leq 1\right).
\eqno(3.2)$$

{\it Proof.} We use the method in the proof of [DZ, Theorem 4.4], which  extends
 the Booss-Wojciechowski theorem to the case of families, to prove (3.2).

Recall that $P_{P_{\partial M}(L)}$ denotes the orthogonal
projection from $L^2(S(TM)\otimes E\otimes {\bf C}^N)$ to the space
of the direct sum of eigenspaces of nonnegative eigenvalues of
$D^E_{P_{\partial M}(L) }$. It is obviously a generalized spectral
section of $D^{\psi, g}_{P_{\partial M}(L)}(u)$ in the sense of \cite{dz}.  Let $P_{P_{\partial M}(L)}(1)$
denote the orthogonal projection from $L^2(S(TM)\otimes E\otimes
{\bf C}^N)$ to the space of the direct sum of eigenspaces of
nonnegative eigenvalues of $D^{\psi,g}_{P_{\partial M}(L) }(1)$.

As in [DZ, (1.11)], let $T(P_{P_{\partial M}(L)}, P_{P_{\partial
M}(L)}(1))$ be the Fredholm operator
$$T\left(P_{P_{\partial M}(L)},
P_{P_{\partial M}(L)}(1)\right)= P_{P_{\partial M}(L)} (1)
P_{P_{\partial M}(L)}:$$
$${\rm Im} \left(
P_{P_{\partial M}(L)}\right)\rightarrow {\rm Im}\left(
P_{P_{\partial M}(L)}(1)\right).\eqno(3.3)$$

Now we observe that the argument in the proof of [DZ, Theorem 4.4]
still works in our present situation, and we obtain,
$$-{\rm sf}\left(D^{\psi,g}_{P_{\partial M}(L)}(u), 0\leq u\leq 1\right)
={\rm ind}\, T\left(P_{P_{\partial M}(L)}, P_{P_{\partial M}(L)}
(1)\right).\eqno(3.4)$$

From (2.7), (3.1), (3.3) and (3.4), one deduces that
$$
-{\rm sf}\left(D^{\psi,g}_{P_{\partial M}(L)}(u), 0\leq u\leq
1\right)
 =  {\rm ind}\, \left(g^{-1}P^{\psi}_{gP_{\partial M}(L) g^{-1}}g
P_{P_{\partial M}(L)}\right)   =  {\rm ind}\, T^E_{g, \psi}(L).
\eqno(3.5)
$$

The proof of Theorem 3.1 is completed.  Q.E.D.

%$\ $

%Theorem 2.1 is our  starting point of computing ${\rm ind}\, T^E_g(L)$. In the next subsection,
%we will reduce the computation of the spectral flow in the right hand side of (2.2) to
%those of variations of $\eta$-invariants.

$\ $

\noindent {\bf b). $\eta$-invariants and the spectral flow}

$\ $

As usual, by \cite{aps1}, for any $u\in [0,1]$, one can define the
$\eta$-invariant ${\eta}(D^{\psi,g}_{P_{\partial M}(L)}(u))$ as well
as the corresponding reduced $\eta$-invariant
$$
\overline{\eta}\left(D^{\psi,g}_{P_{\partial M}(L)}(u)\right)={\dim
\ker\left( D^{\psi,g}_{P_{\partial M}(L)}(u)\right)+\eta
\left(D^{\psi,g}_{P_{\partial M}(L)}(u)\right)\over 2}. \eqno(3.6)
$$

As was mentioned in [KL], it follows from the work of Grubb [Gr]
that when mod ${\bf Z}$, the reduced $\eta$-invariants
$\overline{\eta}(D^{\psi,g}_{P_{\partial M}(L)}(u))$ vary smoothly
with respect to $u\in [0,1]$. And we denote by ${d\over
du}(\overline{\eta}(D^{\psi,g}_{P_{\partial M}(L)}(u))$ the smooth
function on $[0,1]$ of the local variation (after mod ${\bf Z}$)  of
these reduced $\eta$-invariants.

By [KL, Lemma 3.4] and (3.1), one then has
$$
{\rm sf}\left(D^{\psi,g}_{P_{\partial M}(L)}(u), 0\leq u\leq
1\right) =  \overline{\eta}\left(D^{\psi,g}_{P_{\partial
M}(L)}(1)\right)- \overline{\eta}\left(D^{E}_{P_{\partial
M}(L)}\right)   -\int_0^1{d\over
du}\overline{\eta}\left(D^{\psi,g}_{P_{\partial M}(L)}(u)\right)
du.   \eqno(3.7)
$$

By (3.1) and an obvious conjugation, one sees directly that
$$
\overline{\eta}\left(D^{\psi,g}_{P_{\partial M}(L)}(1)\right)
=\overline{\eta}\left(D^{\psi}_{gP_{\partial
M}(L)g^{-1}}\right).\eqno(3.8)
$$

Set $M_{-}=M \setminus ([0, 1]\times \partial M)$. On the boundary
$\partial M_{-}=\{ 1 \}\times
\partial M$ of $M_{-}$, we use the boundary condition
$P_{\partial M}(L)$ and denote it by $P^E_{\partial M_{-}}(L)$. By
[M\"u, Proposition 2.16] one has
$$
\overline{\eta}\left(D^E_{P_{\partial M}(L)}\right)=
\overline{\eta}\left(D^E_{P_{\partial
M_{-}}(L)}\right).\eqno(3.9)
$$

From (2.8), (3.8), (3.9) and using [KL, Theorem 7.7], one deduces
that
$$
\overline{\eta}\left(D^{\psi,g}_{P_{\partial M}(L)}(1)\right)-
\overline{\eta}\left(D^E_{P_{\partial M}(L)}\right) =
\overline{\eta}\left(D^{\psi}_{gP_{\partial M}(L)g^{-1}}\right) -
\overline{\eta}\left(D^{E}_{P_{\partial M_{-}}(L)}\right) \hspace{.7in}
$$
$$ \hspace{1.7in}  =\overline{\eta}\left(D^{\psi,g}_{[0, 1]}\right) - \tau_\mu \left({\cal P}^{\psi}_{[0, 1]},  P_{\partial M}(L) ,
{\cal P}^E_{M_-} \right) ,\eqno(3.10)
$$
where ${\cal P}^E_{M_-}$
is the Calder\'on projection operator associated to $D^E$ on $M_{-}$,
${\cal P}^{\psi}_{[0, 1]}$ the Calder\'on projection operator
associated to $D^{\psi}$ on $[0, 1] \times \partial M$ with the
boundary condition $gP_{\partial M}(L)g^{-1}$ at $\{0\} \times
\partial M$, and
$ \tau_\mu  ({\cal P}^{\psi}_{[0, 1]},  P_{\partial M}(L) , {\cal
P}^E_{M_-}  )
$
is the Maslov triple index in
the sense of Kirk-Lesch [KL, Definition 6.8].

From (3.2), (3.7) and (3.10), one sees that in order to establish an
index formula for $T^E_{g, \psi}$, one needs only to compute
$$
\int_0^1{d\over du}\overline{\eta}\left(D^{\psi,g}_{P_{\partial M}(L)}(u)\right)du.
$$

From (3.1), one verifies that
$${d\over du}D^{\psi,g}_{P_{\partial M}(L)}(u)=(1-\psi)g^{-1}[D^E,g]\eqno(3.12)$$
is a bounded operator.

By the main result in [Gr],  when $t\rightarrow 0^+$, one has the asymptotic
expansion
$${\rm Tr}\left[(1-\psi)g^{-1}[D^E,g]\exp\left(-t D^{\psi,g}_{P_{\partial M}(L)}(u)^2\right)\right]
=\sum_{-\dim M\leq k<0}c_k(u)t^{k/2}$$
$$ +c_0(u)\log t+c_0'(u)+o(1).\eqno(3.13)$$

From (3.12), (3.13) and by proceeding as in [M\"u, Section 2], one
deduces easily that {$$ {d\over
du}\overline{\eta}\left(D^{\psi,g}_{P_{\partial M}(L)}(u)\right)
=-{c_{-1}(u)\over \sqrt{\pi}}.\eqno(3.14)
$$}

The following result gives the explicit value of each $c_{-1}(u)$.

$\ $

\noindent {\bf Theorem 3.2.} {\it We have,
$$
{c_{-1}(u)\over \sqrt{\pi}}=\left({1\over 2\pi\sqrt{-1}}\right)^{(\dim M+1)/2}\int_M
\widehat{A}\left(R^{TM}\right){\rm Tr}\left[\exp\left(-R^E\right)\right]$$
$$\cdot \, {\rm Tr}\left[g^{-1}dg\exp\left( (1-u)u \left(g^{-1}dg\right)^2\right)
\right].\eqno(3.15)
$$
Therefore,
$$
 {\rm ind}\, T^E_{g, \psi} = -
\left({1\over 2\pi\sqrt{-1}}\right)^{(\dim M+1)/2}\int_M
\widehat{A}\left(R^{TM}\right){\rm
Tr}\left[\exp\left(-R^E\right)\right]{\rm ch}(g)
$$
$$
-  \overline{\eta}\left(D^{\psi,g}_{ {[0, 1]}}\right) + \tau_\mu
\left({\cal P}^{\psi}_{[0, 1]},  P_{\partial M}(L) , {\cal
P}^E_{M_-} \right) . \hspace{.5in}  \eqno(3.16)
$$}

%$\ $

\noindent {\bf Remark 3.3.} When $g|_{\partial M}={\rm Id}$,
Theorem 3.2 was proved by Douglas and Wojciechowski \cite{dw}.
Thus, our main concern will be the case where $g|_{\partial M}$ is
not identity. Note that in this case, $gP_{\partial
M}(L)g^{-1}-P_{\partial M}(L)$ is in general not a smoothing
operator.

$\ $

\noindent {\bf c). A Proof of Theorem 3.2}

$\ $

Here we prove Theorem 3.2. The contribution to the left hand side of
(3.15) splits into the interior and boundary parts. The interior
contribution can be handled by the standard local index theory
techniques and the boundary contribution can be easily seen to be
zero.

Recall that our purpose is to study the asymptotic behavior when
$t\rightarrow 0^+$ of the following quantity:
$$
{\rm Tr} \left[(1-\psi) g^{-1}\left[ D^E,g\right]\exp\left(-t
\left(D^{\psi,g}_{P_{\partial M}(L)}(u)\right)^2\right)\right],\eqno(3.17)
$$
for  $0\leq u\leq 1$,
where $D^{\psi, g}_{P_{\partial M}(L)}(u)$ is defined by (3.1):
$$
D_{\psi}(u)=(1-u)D^E+ug^{-1}D^{\psi}g=D^E+u(1-\psi)g^{-1}\left[
D^E,g\right]. \eqno(3.18)
$$

From (3.18) one verifies that
$$
\left( D_{\psi}(u)\right)^2=\left(D^E\right)^2+B(u),\eqno(3.19)
$$
where
$$
B(u)=u\left[D^E,(1-\psi)g^{-1}\left[ D^E,g\right]\right]
+u^2(1-\psi)^2\left(g^{-1}\left[ D^E,g\right]\right)^2.\eqno(3.20)
$$
Since
$$
\left[D^E,g^{-1}\left[ D^E,g\right]\right] =\left[ D^E,g^{-1}\right]\left[ D^E,g\right]
+g^{-1}\left[D^E,\left[ D^E,g\right]\right]
$$
$$
\hspace{1.2in}  =-\left(g^{-1}\left[ D^E,g\right]\right)^2+g^{-1}\left[D^E,\left[ D^E,g\right]\right],\eqno(3.21)
$$
one finds that
$$
B(u)=\left(u^2(1-\psi)^2-u(1-\psi)\right)\left(g^{-1}\left[
D^E,g\right]\right)^2 + u(1-\psi)g^{-1}\left[D^E,\left[
D^E,g\right]\right]
$$
$$ - u\psi' c(\frac{\partial}{\partial x})
g^{-1}\left[ D^E,g\right]. \hspace{2.51in}  \eqno(3.22)
$$

Let $e_1,\dots,e_{\dim M}$ be an orthonormal basis of $TM$. Then by (2.3), one verifies that
$$
\left[D^E,\left[ D^E,g\right]\right]=\sum_{i,j=1}^{\dim M}
\left[ c(e_i)\nabla_{e_i},c(e_j)\left(\nabla_{e_j}g\right)\right] \hspace{2.1in}
$$
$$
=-\sum_{i=1}^{\dim M} \left(\nabla_{e_i}^2-\nabla_{\nabla^{TM}_{e_i}e_i}\right)g-2\sum_{i=1}^{\dim M}
\left(\nabla_{e_i}g\right)\nabla_{e_i}.\eqno(3.23)
$$ Therefore,
$B(u)$ is a differential operator of order one.

For brevity of notation, from now on we will denote the elliptic
operator $D^{\psi,g}_{P_{\partial M}(L)}(u)$ simply by $D_{\psi}(u)$
(with the boundary condition $P_{\partial M}(L)$ understood).

Also, denote $\dim M=2n+1$. We first show that the study of the
limit as $t\rightarrow 0^+$ of the term in (3.15) can be reduced to
separate computations in the interior and near the boundary. We fix the $\epsilon$
which defines the cut off function $\psi$ (for example, we can take $\epsilon = \frac{1}{4}$).

Let $\overline{D}_{\psi}(u)$ be the double of the Dirac type
operator $D_{\psi}(u)$, which lives on the double of $M$. Let
$E_I(t)$ denote the heat kernel associated to
$e^{-t(\overline{D}_{\psi}(u))^2}$. Let $E_{L,b}(t)$ denote the heat
kernel of $e^{-t(D^{\psi}_{gP_{\partial M}(L)g^{-1}})^2}$ on the
half cylinder $[0,+\infty)\times \partial M$, where we extend
everything from $[0,\epsilon]\times \partial M$  canonically. By our
assumption, this is simply the conjugation of the heat kernel
$e^{-t(D^{E}_{P_{\partial M}(L)})^2}$ on the half cylinder:
$$
e^{-t(D^{\psi}_{gP_{\partial M}(L)g^{-1}})^2}= g
e^{-t(D^{E}_{P_{\partial M}(L)})^2} g^{-1} , \eqno(3.24)
$$
where $D^E$ assumes the product form (2.4).

Following  \cite{aps1}, we use $\rho(a, b)$ to denote an increasing $C^{\infty}$ function
of the real variable $x$ such that
$$\rho=0 \ \ \mbox{for} \ x\leq a, \ \  \ \rho=1 \ \ \mbox{for} \ x\geq b.  \eqno(3.25)$$
Define four $C^{\infty}$ functions by
$$
 \phi_1 = 1 -\rho\left(\o{5}{6}\epsilon, \epsilon \right), \ \  \psi_1 = 1 -\rho\left(\o{3}{6}\epsilon, \o{4}{6} \epsilon \right),
 $$
 $$\phi_2 = \rho\left(\o{1}{6}\epsilon,  \o{2}{6}\epsilon \right), \ \  \psi_2 = \rho\left(\o{3}{6}\epsilon, \o{4}{6}\epsilon  \right)
 . \eqno(3.26)
$$

\noindent {\bf Lemma 3.4.} {\it There exists $C>0$ such that as $t\rightarrow 0^+$, one has,}
$$
{\rm Tr}\left[(1-\psi)g^{-1}
 \left[D^E,g\right] \exp\left(-t\left(D_{\psi}(u)\right)^2\right)\right] \hspace{2.5in}
 $$
$$ =\int_M{\rm Tr}\left[(1-\psi) g^{-1}
 \left[D^E,g\right]E_I(t)(x,x)\right]\psi_2(x)d{\rm vol} \hspace{1.8in}
 $$
 $$
\hspace{.5in} +\int_{[0,+\infty)\times \partial M}{\rm Tr}\left[(1-\psi) g^{-1}
 \left[D^E,g\right]E_{L,b}(t)(x,x)\right]\psi_1(x)d{\rm vol}
 +O\left(e^{-C/t}\right).\eqno(3.27)
 $$

 {\it Proof.} We construct a parametrix for $\exp\left(-t\left(D_{\psi}(u)\right)^2\right)$ by patching:
$$
E(t)= \phi_1 E_{L,b}(t) \psi_1 + \phi_2 E_I(t) \psi_2.\eqno(3.28)
$$
By the standard theory the interior heat kernel is exponentially small as $t \ra 0^+$
for $x\not= y$. That is, there exists $C_1>0$ such that for $0<t \leq 1$ (say),
$$
|E_I(t)(x, y)| \leq C_1 t^{-n-\o{1}{2}} e^{-\o{d(x,y)^2}{4t}}. \eqno(3.28)
$$
Moreover the same estimate holds for derivatives of $E_I(t)(x,y)$ if we replace
$ t^{-n-\o{1}{2}}$ by $ t^{-n-\o{1}{2}-l_1-\o{l_2}{2}}$ where $l_1$ is the number of
time differentiation and $l_2$ is the number of spatial differentiation.

One has the same estimate for $E_{L, b}(t)$ as shown in \cite[Proposition 2.21]{aps1}):
$$
|E_{L,b}(t)(x, y)| \leq C_2 t^{-n-\o{1}{2}} e^{-\o{d(x,y)^2}{4t}}.\eqno(3.29)
$$
Furthermore,  similar estimates for derivatives continue to hold as in the case of interior heat kernel.

By our construction, the distance of the support of $\phi_i'$, $i=1,\ 2$, to the support of
$\psi_i$, $i=1,\ 2$, is at least ${1\over 6} \epsilon$. Therefore the estimates above give
$$
 \left(\o{\pt}{\pt t} + \left(D_{\psi}(u)\right)^2\right)E(t)= O\left(e^{-\o{C}{t}}\right),\eqno(3.30)
 $$
for some $C=C(\epsilon)>0$, and the derivatives of $(\o{\pt}{\pt t}
+ (D_{\psi}(u))^2)E(t)$ decays exponentially as well (with a smaller
$C$). Hence by Duhamel principle, one deduces that
\begin{align*}
\left[(1-\psi) g^{-1}
 \left[D^E,g\right] \exp\left(-t\left(D_{\psi}(u)\right)^2\right)\right]
 = & \phi_1 \left[(1-\psi) g^{-1} \left[D^E,g\right] E_{L,b}(t)\right]
 \psi_1 \\
& + \phi_2 \left[(1-\psi) g^{-1} \left[D^E,g\right] E_I(t)\right]
\psi_1 + O\left(e^{-\o{C}{t}}\right).
 \end{align*}
Our result follows.  Q.E.D.

 $\ $

 Clearly,
$$
{\rm Tr}\left[(1-\psi) g^{-1}
 \left[D^E,g\right]E_{L,b}(t)(x,x)\right]\psi_2(x)=0.\eqno(3.31)
$$
We therefore turn our attention to the interior contribution.

$\ $

\noindent {\bf Lemma 3.5.} {\it We have, $t\rightarrow 0^+$,}
$$
\int_M {\rm Tr}\left[(1-\psi) g^{-1}
 \left[D^E,g\right]E_I(t)(x,x)\right]\psi_2 d{\rm vol} \rightarrow
 $$
 $$
 \left({1\over 2\pi\sqrt{-1}}\right)^{n+1}\int_M
\widehat{A}\left(R^{TM}\right){\rm Tr}\left[
\exp\left(-R^E\right)\right]{\rm ch}(g)d{\rm vol}. \eqno(3.32)
$$

{\it Proof.} By (3.22) and by applying by now the standard local
index techniques analogous to [G] and [DZ, Section 4e)], we obtain
that as $t\rightarrow 0^+$,
 $$
 {\rm Tr}\left[(1-\psi) g^{-1}
 \left[D^E,g\right]E_I(t)(x,x)\right]  \rightarrow
  \left({1\over 2\pi\sqrt{-1}}\right)^{n+1}
\widehat{A}\left(R^{TM}\right){\rm Tr}\left[
\exp\left(-R^E\right)\right] \cdot
$$
$$
\int_0^1 {\rm Tr} \left[ (1-\psi) g^{-1}dg \exp
\left(\left(u(1-\psi)-u^2(1-\psi)^2\right) (g^{-1}dg)^2 + u d\psi
g^{-1}dg \right)\right]du  .\eqno(3.33)
$$
It follows from the nilpotency of $d\psi$ that
$$
{\rm Tr} \left[ (1-\psi) g^{-1}dg \exp
\left(\left(u(1-\psi)-u^2(1-\psi)^2\right)(g^{-1}dg)^2 + u d\psi
g^{-1}dg \right)\right] \hspace{.2in}
$$
$$
= {\rm Tr} \left[ (1-\psi) g^{-1}dg \exp
\left(\left(u(1-\psi)-u^2(1-\psi)^2\right) (g^{-1}dg)^2
\right)\right] \hspace{1in}
$$
$$
+ {\rm Tr} \left[ (1-\psi)u d\psi \left(g^{-1}dg\right)^2 \exp
\left(\left(u(1-\psi)-u^2(1-\psi)^2\right) (g^{-1}dg)^2
\right)\right]. \eqno(3.34)
$$

Since
$$
{\rm Tr} \left[ \left(g^{-1}dg\right)^{2k} \right]=0 ,\eqno(3.35)
$$
for any positive integer $k$, the second term on the right hand side
of (3.34) is zero. On the other hand, the first term on the right
hand side of (3.34), when restricted to the cylindrical part of $M$,
contains no form in the normal direction $x$ by our product
structure assumption on $g$. Thus its integration over the
cylindrical part of $M$, where the cut off function $\psi$ may not
be zero, is zero. This is true even when integrated together with
$\widehat{A}\left(R^{TM}\right){\rm Tr}\left[
\exp\left(-R^E\right)\right]$ as we also have product structure
assumption on these geometric data. It follows then that the right
hand side of (3.33) equals
$$
 \left({1\over 2\pi\sqrt{-1}}\right)^{n+1}\int_M
\widehat{A}\left(R^{TM}\right){\rm Tr}\left[
\exp\left(-R^E\right)\right]   \int_0^1 {\rm Tr} \left[ g^{-1}dg
\exp \left((1-u)u (g^{-1}dg)^2 \right)\right]du .\eqno(3.36)
$$
Lemma 3.5 follows. Q.E.D.

By Lemmas 3.4, 3.5, one gets (3.15). Then (3.16) follows from
(3.4), (3.5), (3.7), (3.10), (3.14) and (3.15).

The proof of Theorem 3.2 is now complete.\ \ Q.E.D.

\section{Spectral flow, Maslov indices and the index of the Toeplitz operator}

In this section, we prove the index formula for the Toeplitz
operator $T^E_g$, as stated in Theorem 2.3. As we have seen in the
previous section, an index formula (3.16) for the perturbed Toeplitz
operator $T^E_{g, \psi}$ has been established. To go from the
perturbed Toeplitz operator to the original Toeplitz operator, we
make use of the spectral flow, reformulated in \cite{dz} in terms of
generalized spectral section, and the theory of Maslov indices, as
developed by \cite{kl}.

$\ $

\noindent {\bf a). Comparison of indices of Toeplitz and perturbed Toeplitz operators}

$\ $

Here we show that the difference of the index of the Toeplitz operator and that of the perturbed Toeplitz operator
can be expressed in terms of a spectral flow by using the formulation of \cite{dz} via generalized spectral sections.

$\ $

 \noindent {\bf Lemma 4.1.} {\it We have,
$$
{\rm ind}\, T^E_g - {\rm ind}\, T^E_{g, \psi} =  {\rm sf} \left(
D^{\psi}(s), 0\leq s \leq 1 \right), \eqno(4.1)
$$ where
$$
D^{\psi}(s)= g D^{\psi,g}(s) g^{-1}
=g\left(D^E+(1-s\psi)g^{-1}[D^E,g]\right)g^{-1} \eqno(4.2)
$$
is   equipped with the boundary condition $g P_{\partial M}(L)
g^{-1}$. }

 {\it Proof.} First, we note that
$$ {\rm ind}\, T^E_g(L)={\rm ind} \left( P_{gP_{\partial M}(L) g^{-1}}
g P_{P_{\partial M}(L)} \right) = {\rm ind} \left( P^g_{P_{\partial
M}(L)} P_{P_{\partial M}(L)} \right)$$ where $P^g_{P_{\partial
M}(L)}$ is the orthogonal projection onto the eigenspaces of
$(g^{-1} D^E g, P_{\partial M})$ with nonnegative eigenvalues.
Similarly,
$$ {\rm ind}\, T^E_{g, \psi}(L)={\rm ind} \left( P^{\psi}_{gP_{\partial M}(L) g^{-1}}
g P_{P_{\partial M}(L)} \right) = {\rm ind} \left( P^{g,
\psi}_{P_{\partial M}(L)} P_{P_{\partial M}(L)} \right)$$ where
$P^{g, \psi}_{P_{\partial M}(L)}$ is the orthogonal projection
onto the eigenspaces of $(g^{-1} D^{\psi} g, P_{\partial M})$ with
nonnegative eigenvalues.

Thus,
$${\rm ind}\, T^E_g - {\rm ind}\, T^E_{g, \psi} = {\rm ind} \left( P^g_{P_{\partial
M}(L)} P_{P_{\partial M}(L)} \right) - {\rm ind} \left( P^{g,
\psi}_{P_{\partial M}(L)} P_{P_{\partial M}(L)} \right). $$

Noting that $P_{P_{\partial M}(L)}$ is again a generalized
spectral section of $D^{\psi,g}(s)$, we have by the argument in
\cite{dz} that
$$
{\rm ind}\, T^E_g - {\rm ind}\, T^E_{g, \psi} =  {\rm sf} \left(
D^{\psi,g}(s), 0\leq s \leq 1 \right) =  {\rm sf} \left(
D^{\psi}(s), 0\leq s \leq 1 \right).
$$ Q.E.D.

$\ $

\noindent {\bf b). Maslov indices and the splitting of spectral
flow}

$\ $

Already from Theorem 3.2 and Lemma 4.1 we obtain an index formula for
$T^E_g$.  To put this formula into the (much better) form as stated in our main result,
Theorem 2.3, we need to make use of Maslov indices as developed in \cite{kl}.

The (double) Maslov index is an integer invariant for Fredholm pairs
of paths of Lagrangian subspaces. It is an algebraic count of how
many times these Lagrangian subspaces intersect along the path. We
will follow the treatment of \cite{kl} closely.

Let $H$ be an Hermitian symplectic Hilbert space, i.e., there is an
unitary map $J: \ H \rightarrow H$ such that $J^2=-1$ and the
eigenspaces with eigenvalues $\pm \sqrt{-1}$ have equal dimension.
The Lagrangian subspaces in $H$ can be identified with their
orthogonal projections, the space of which is
$$ {\rm Gr}(H)=\left\{ P \in B(H) \ | \ P=P^*, \ P^2 = P, \ JPJ^* =
I -P\right\} . $$ A pair $(P, Q)$, $P, Q \in {\rm Gr}(H)$ is
called Fredholm if
$$ T(Q,P)=PQ: \ {\rm Im}\, Q \rightarrow {\rm Im}\, P $$
is Fredholm.

For a (continuous) path $(P(t), Q(t))$, $0\leq t \leq 1$,  of
Fredholm pairs, $P(t), Q(t) \in {\rm Gr}(H)$, the Maslov index
associates an integer ${\rm Mas}(P(t), Q(t))$ \cite{kl}.

On the other hand, for a triple $P, Q, R \in {\rm Gr}(H)$ such that
$(P, Q)$, $(Q, R)$, $(P, R)$ are Fredholm and at least one of the
differences $P-Q$, $Q-R$, $P-R$ is compact, an integer
$\tau_{\mu}(P, Q, R)$ can be defined \cite{kl}, which is called the
Maslov triple index. They satisfy the following important relation
\cite[(6.24)]{kl}\footnote{Note the sign correction on the left hand
side of (4.3)}
$$
 \tau_{\mu}(P(1), Q(1), R(1))-  \tau_{\mu}(P(0), Q(0), R(0))
$$
$$
={\rm Mas}(P(t), Q(t)) + {\rm Mas}(Q(t), R(t)) - {\rm Mas}(P(t),
R(t)).  \eqno(4.3)
$$

The following theorem is a slight generalization of \cite[Theorem
7.6]{kl}, which itself is a generalization of a result of
Nicolaescu. It follows from the same argument.

$\ $

\noindent {\bf Theorem 4.2.} {\it Let $M$ be a manifold with
boundary and $H$ a separating hypersurface in $M$ such that $H\cap
\partial M=\emptyset$ and that $M=M^+ \cup_H M^-$. Let $D(t), \ 0\leq t \leq 1$, be a smooth
of Dirac type operators equipped with self adjoint elliptic
boundary conditions of APS type on the boundary of $M$. If $D(t)$
is of product type near the separating hypersurface $H$, then
$$
{\rm sf} \left( D(t), \ 0\leq t \leq 1 \right) = \Mas \left( {\cal
P}_{M^-}, {\cal P}_{M^+} \right), \eqno(4.4)
$$
 where ${\cal P}_{M^-}$ (${\cal
P}_{M^+}$) denotes the Calder\'on projections on $M^-$ ($M^+$) with
the boundary conditions on $M^- \cap \bM$ ($M^+ \cap \bM$) coming
from those on $\bM$.}

$\ $

\noindent {\bf c). Proof of Theorem 2.3}

$\ $

We are now in position to prove our main result. By Lemma 4.1 we have
\[ \ind\, T^E_g - {\rm ind}\, T^E_{g, \psi} =  {\rm sf} \left( D^{\psi}(s),
0\leq s \leq 1 \right). \] Applying Theorem 4.2 to $M=[0, 1] \times
\bM \cup M_-$ with the boundary condition $g P_{\bM}(L) g^{-1}$ on
$\bM$, we obtain
$$
 {\rm sf} \left( D^{\psi}(s), 0\leq s \leq 1 \right)= \Mas \left( {\cal
P}^{\psi}_{[0, 1]}(s), {\cal P}_{M_-} \right). \eqno(4.5)
$$
Here ${\cal P}^{\psi}_{[0, 1]}(s)$ denotes the Calder\'on projection
operator associated to $D^{\psi}(s)$ on $[0, 1] \times \partial M$
with the boundary condition $gP_{\partial M}(L)g^{-1}$ at $\{0\}
\times \partial M$

Hence by Theorem 3.2 and (4.5), we have
$$
\ind\, T^E_g=  - \left({1\over 2\pi\sqrt{-1}}\right)^{(\dim
M+1)/2}\int_M \widehat{A}\left(R^{TM}\right){\rm
Tr}\left[\exp\left(-R^E\right)\right]{\rm ch}(g) \hspace{.8in}
$$
$$
- \overline{\eta}\left(D^{\psi,g}_{[0, 1]}\right) +   \tau_\mu
\left({\cal P}^{\psi}_{[0, 1]},  P_{\partial M}(L) , {\cal
P}_{M_-} \right)  + \Mas \left( {\cal P}^{\psi}_{[0, 1]}(s),
{\cal P}_{M_-} \right). \eqno(4.6)
$$
Using (2.13), we rewrite (4.6) as
$$
\ind\, T^E_g=  - \left({1\over 2\pi\sqrt{-1}}\right)^{(\dim
M+1)/2}\int_M \widehat{A}\left(R^{TM}\right){\rm
Tr}\left[\exp\left(-R^E\right)\right]{\rm ch}(g) -
\overline{\eta}(\bM, g) \hspace{.8in}
$$
$$
\hspace{.6in}  -{\rm sf} \{D^{\psi,g}_{[0, 1]}(s); 0 \leq s \leq 1
\} +
 \tau_\mu \left({\cal P}^{\psi}_{[0, 1]},  P_{\partial M}(L) ,
{\cal P}_{M_-} \right)  + \Mas \left( {\cal P}^{\psi}_{[0, 1]}(s),
{\cal P}_{M_-} \right). \eqno(4.7)
$$

On the other hand, by (4.3),
$$
 \tau_{\mu} ( {\cal P}^{\psi}_{[0, 1]},  P_{\bM}(L), {\cal P}_{M_-})
 - \tau_{\mu} (g P_{\bM}(L) g^{-1}, P_{\bM}(L), {\cal P}_{M_-}) =
$$
$$
 \Mas ( {\cal P}^{\psi}_{[0, 1]}(s),  P_{\bM}(L)) - \Mas (
{\cal P}^{\psi}_{[0, 1]}(s), {\cal P}_{M_-}) . \eqno(4.8)
$$
And finally, by using \cite[Theorem 7.5]{kl},
$$
\hspace{.5in}  {\rm sf} \{D^{\psi,g}_{[0, 1]}(s); 0 \leq s \leq 1
\}={\rm sf} \{D^{\psi}_{[0, 1]}(s); 0 \leq s \leq 1 \}=\Mas (
{\cal P}^{\psi}_{[0, 1]}(s), P_{\bM}(L)). \eqno(4.9)
$$

From (4.7)-(4.9), one gets
$$
\ind\, T^E_g=  - \left({1\over 2\pi\sqrt{-1}}\right)^{(\dim
M+1)/2}\int_M \widehat{A}\left(R^{TM}\right){\rm
Tr}\left[\exp\left(-R^E\right)\right]{\rm ch}(g) -
\overline{\eta}(\bM, g) \hspace{.8in}
$$
$$
  + \tau_{\mu} (g P_{\bM}(L) g^{-1}, P_{\bM}(L), {\cal
P}_{M_-}). \eqno(4.10)
$$
This completes the proof of Theorem 2.3.\ \
Q.E.D.

\section{Generalizations and some further results}

In this section, we first show that for any even dimensional
closed spin manifold $X$ and any $K^1$ representative $g:\ X \ra
U(N)$, the invariant $\overline{\eta}(X, g)$ defined in (2.13) is
independent of the cut off function. Then we generalize Theorem
2.3 to the case where one no longer assumes that $g$ is of the
product type near $\partial M$. Finally we take a further look at
the $\eta$-type invariant $\overline{\eta}(X, g)$ and study some
of its basic properties.

This section is organized as follows. In a), we study the variation
of  $\overline{\eta}(X, g)$ in the cut off function which gives us
the desired independence. We also make a conjecture about what this
eta invariant really is. In b), we take a look at the variations of
the odd Chern character  forms . In c), we prove an extension of
Theorem 1.3 to the case where we no longer assume $g$ is of product
structure near $\partial M$. In d), we make a further study of the
eta invariant $\overline{\eta}(X, g)$.

$\ $

\noindent {\bf a). The invariant $\overline{\eta}(X, g)$}

$\ $

Recall that the invariant of $\eta$ type associated to a Dirac
operator on an   even dimensional manifold $X$ with vanishing
index and the $K^1$ representative $g$ over $X$ is defined in
(2.13) as (here we have inserted $\psi$ in the notation to
indicate that, a priori, it depends on the cut off function
$\psi$)
$$
\overline{\eta}(X, g, \psi)=
 \overline{\eta}(D^{\psi, g}_{[0, 1]}) -
 {\rm sf} \left\{D^{\psi,g}_{[0, 1]}(s); 0 \leq s \leq 1 \right\},
$$
where $D^{\psi,g}_{[0, 1]}(s)$ is a path connecting $g^{-1} D^E g$
with $D^{\psi,g}_{[0, 1]}$ defined by
$$
D^{\psi,g}(s)= D^E + (1-s\psi) g^{-1}[D^E, g]
$$
on $[0, 1]\times X$, with the boundary condition $P_{X}(L) $ on
$\{0\}\times X$ and the boundary condition ${\rm Id}-g^{-1}P_{X}(L)
g$ at $\{1\}\times X$.

$\ $

\noindent {\bf Proposition 5.1.} {\it The invariant
$\overline{\eta}(X, g, \psi)$ is independent of the cut off function
$\psi$.}

$\ $

{\it Proof.} Let $\psi_1$, $\psi_2$ be two cut off functions and
$$\psi_t=(2-t)\psi_1+(t-1)\psi_2,\ \ \ \  \ 1\leq t \leq 2,\eqno(5.1)$$ be the smooth path of cut off functions
connecting the two. Then
$$
\overline{\eta}(D^{\psi_2, g}_{[0, 1]}) -
\overline{\eta}(D^{\psi_1, g}_{[0, 1]})= \int_1^2
{\partial\over\partial t} \overline{\eta}(D^{\psi_t, g}_{[0, 1]})
dt + {\rm sf} \left\{ D^{\psi_t, g}_{[0, 1]}, 1\leq t \leq 2
\right\}. \eqno(5.2)
$$
As before, we can compute ${\partial\over\partial t}
\overline{\eta}(D^{\psi_t, g}_{[0, 1]})$ via heat kernel and local
index theorem technique (Cf. Section 3 {\bf c)}) and find
$$
{\partial\over\partial t} \overline{\eta}(D^{\psi_t, g}_{[0, 1]})
\equiv 0. \eqno(5.3)
$$
Here we have once again used the fact that $g$ is constantly
extended along the radial direction. Therefore
$$
\overline{\eta}(X, g, \psi_2) - \overline{\eta}(X, g, \psi_1)=
{\rm sf} \left\{ D^{\psi_t, g}_{[0, 1]}, 1\leq t \leq 2 \right\} -
{\rm sf} \left\{D^{\psi_2,g}_{[0, 1]}(s); 0 \leq s \leq 1 \right\}
$$
$$
+{\rm sf} \left\{D^{\psi_1,g}_{[0, 1]}(s); 0 \leq s \leq 1
\right\} =0
$$
by the additivity and the homotopy invariance of spectral flow.
  Q.E.D.

$\ $

 Thus, the eta type invariant $\overline{\eta}(X, g)$, which
we introduced using a cut off function, is in fact independent of
the cut off function. This leads naturally to the question of
whether $\overline{\eta}(X, g)$ can actually be defined directly.
We now state a conjecture for this question.

Let $D^{[0, 1]}$ be the Dirac operator on $[0, 1]\times X$. We
equip the boundary condition $gP_{X}(L) g^{-1}$ at $\{0\}\times X$
and the boundary condition ${\rm Id}-P_{X}(L)$ at $\{ 1 \}\times
X$.

Then $(D^{[0, 1]}, gP_{X}(L) g^{-1} , {\rm Id}-P_{X}(L)
 )$ forms a self-adjoint elliptic boundary problem. We denote
the corresponding elliptic self-adjoint operator by
$D^{[0,1]}_{gP_{X}(L) g^{-1} , P_{X}(L) }$.

Let $\eta(D^{[0,1]}_{gP_{X}(L) g^{-1}, P_{X}(L)  },s)$ be the
$\eta$-function of $D^{[0,1]}_{gP_{X}(L) g^{-1} , P_{X}(L)
 }$. By [KL, Theorem 3.1], which goes back to [Gr], one knows
that the $\eta$-function $\eta(D^{[0,1]}_{gP_{X}(L) g^{-1} ,
P_{X}(L)  },s)$ admits a meromorphic extension to ${\bf C}$ with
poles of order at most 2. One then defines, as in [KL, Definition
3.2], the $\eta$-invariant of $D^{[0,1]}_{gP_{X}(L) g^{-1},
P_{X}(L)  }$, denoted by $\eta(D^{[0,1]}_{gP_{X}(L) g^{-1} ,
P_{X}(L)  })$, to be the constant term in the Laurent expansion of
$\eta(D^{[0,1]}_{gP_{X}(L) g^{-1}, P_{X}(L)  },s)$ at $s=0$.

Let $\overline{\eta}(D^{[0,1]}_{gP_{X}(L) g^{-1}, P_{X}(L)
 })$ be the associated reduced $\eta$-invariant.

$\ $

\noindent {\bf Conjecture 5.2}:
$$
\overline{\eta}(X, g)= \overline{\eta}(D^{[0,1]}_{gP_{X}(L) g^{-1}
, P_{X}(L)  }).
$$

If this conjecture is correct, then the result stated in [Z2,
Theorem 5.2] is also correct. A previous version of the current
article was devoted to a proof of [Z2, Theorem 5.2], and a referee
pointed out a gap in that version. This is why we now introduce a
new $\eta$-type invariant, which makes the   picture clearer.

 $\ $

\noindent {\bf  b). A Chern-Weil type theorem for odd Chern character forms}

$\ $

In this subsection, we assume that there is a smooth  family $g_t$, $0\leq t\leq 1$, of
the automorphisms of the trivial complex vector bundle ${\bf C}^N\rightarrow M$, and study
the variations of the odd Chern character forms ${\rm ch}(g_t,d)$, when $t\in [0,1]$
changes.

The following lemma is taken from [G, Proposition 1.3] (cf. [Z1,
Lemma 1.17]).

$\ $

\noindent {\bf Lemma 5.3.} {\it For any positive odd integer $n$,
the following identity holds,}
$${\partial\over\partial t}{\rm Tr}\left[\left(g_t^{-1}dg_t\right)^n\right]=
nd{\rm Tr}\left[g_t^{-1}{\partial g_t\over \partial t}
\left(g_t^{-1}dg_t\right)^{n-1}\right].\eqno(5.4)$$

{\it Proof.} %We include a proof for the self-contain of the present paper.
First of all, from the identity $g_tg_t^{-1}={\rm Id}$, one verifies
by differentiation that
$${\partial g_t^{-1}\over\partial t}=-g_t^{-1}\left(
{\partial g_t\over \partial t}\right)g_t^{-1}.\eqno(5.5)$$

One then computes that
$${\partial\over\partial t}\left(g_t^{-1}dg_t\right)=
{\partial g_t^{-1}\over\partial t}dg_t+g_t^{-1}d{\partial g_t\over\partial t}
=-\left(g_t^{-1}{\partial g_t\over\partial t}\right)g_t^{-1}dg_t+
g_t^{-1}d{\partial g_t\over\partial t}$$
$$=-\left(g_t^{-1}{\partial g_t\over\partial t}\right)g_t^{-1}dg_t
+\left(g_t^{-1}dg_t\right)\left(g_t^{-1}{\partial g_t\over\partial
t}\right) +d\left(g_t^{-1}{\partial g_t\over\partial
t}\right).\eqno(5.6)$$

One also verifies that
$$d\left(g_t^{-1}dg_t\right)^2=d\left(g_t^{-1}dg_t\right)g_t^{-1}dg_t
-g_t^{-1}dg_td\left(g_t^{-1}dg_t\right)
=0,$$
from which one deduces that for any positive even integer $k$,
$$d\left(g_t^{-1}dg_t\right)^k=0.\eqno(5.7)$$

From (5.5)-(5.7), one verifies that
$${\partial\over\partial t}{\rm Tr}\left[\left(g_t^{-1}dg_t\right)^n\right]
=n{\rm Tr}\left[
{\partial\over\partial t}\left(g_t^{-1}dg_t\right)\left(g_t^{-1}dg_t\right)^{n-1}\right]$$
$$=n{\rm Tr}\left[\left[g_t^{-1}dg_t,g_t^{-1}{\partial g_t\over\partial t}\right]
\left(g_t^{-1}dg_t\right)^{n-1}\right]
+n{\rm Tr}\left[
d\left(g_t^{-1}{\partial g_t\over\partial t}\right)\left(g_t^{-1}dg_t\right)^{n-1}\right]$$
$$=n{\rm Tr}\left[\left[g_t^{-1}dg_t,g_t^{-1}{\partial g_t\over\partial t}
\left(g_t^{-1}dg_t\right)^{n-1}\right]\right]
+n{\rm Tr}\left[
d\left(g_t^{-1}{\partial g_t\over\partial t}\left(g_t^{-1}dg_t\right)^{n-1}
\right)\right]$$
$$=nd{\rm Tr}\left[g_t^{-1}{\partial g_t\over\partial t}\left(g_t^{-1}dg_t\right)^{n-1}
\right].\eqno(5.8)$$

The proof of Lemma 5.3 is completed.  Q.E.D.

$\ $

For any $t\in [0,1]$, set
$$\widetilde{\rm ch}\left(g_t,d,0\leq t\leq 1\right)=
\sum_{n=0}^{(\dim M-1)/2}{n!\over (2n)!}\int_0^1{\rm
Tr}\left[g_t^{-1} {\partial g_t\over\partial
t}\left(g_t^{-1}dg_t\right)^{2n}\right]dt.\eqno(5.9)$$

By Lemma 5.3  and (5.9), one gets

$\ $

\noindent {\bf Theorem 5.4.} (cf. [G, Proposition 1.3]) {\it The
following identity holds,}
$${\rm ch}(g_1,d)-{\rm ch}(g_0,d)=d\widetilde{\rm ch}\left(g_t,d,0\leq t\leq 1\right).\eqno
(5.10)$$

$\ $

\noindent {\bf c). An index theorem for the case of non-product
structure near boundary}

$\ $

In this section, we no longer assume the product structure of $g:M\rightarrow U(N)$ near
the boundary $\partial M$. Then, clearly, the Toeplitz operator $T_g^E(L)$ is still
well-defined. Moreover, by an easy deformation argument, we can construct a smooth one
parameter family of maps $g_t:M\rightarrow U(N)$, $0\leq t\leq 1$, with $g_0=g$,
$g_1=g'$ such that for any $t\in [0,1]$, $g_t|_{\partial M}=g_0|_{\partial M}$, and
that $g_1=g'$ is of product structure near $\partial M$.

By the homotopy invariance of the index of Fredholm operators, one has
$${\rm ind}\, T_g^E(L) ={\rm ind}\, T_{g'}^E(L) .\eqno(5.11)$$

Now by Theorem 2.3, one has
$${\rm ind}\, T_{g'}^E(L) = -\left({1\over 2\pi\sqrt{-1}}\right)^{(\dim M+1)/2}\int_M
\widehat{A}\left(R^{TM}\right){\rm Tr}\left[\exp\left(-R^E\right)\right]{\rm ch}(g',d)$$
$$- \eta\left(\bM, g\right)+ \tau_\mu \left(gP_{\partial M}(L)g^{-1} , P_{\partial
M}(L)  , {\cal P}_M\right) .\eqno(5.12)$$

Since $g_t|_{\partial M}$ is constant in $t$, from (5.9) and (5.10)
one deduces that
$$\int_M
\widehat{A}\left(R^{TM}\right){\rm Tr}\left[\exp\left(-R^E\right)\right]{\rm ch}(g',d)
-\int_M
\widehat{A}\left(R^{TM}\right){\rm Tr}\left[\exp\left(-R^E\right)\right]{\rm ch}(g,d)$$
$$=\int_{\partial M}
\widehat{A}\left(R^{TM}\right){\rm
Tr}\left[\exp\left(-R^E\right)\right] \widetilde{\rm ch}(g_t,d,0\leq
t\leq 1)=0.\eqno(5.13)$$

From (5.11)-(5.13), one deduces

$\ $

\noindent{\bf Theorem 5.5.} {\it Formula (2.16) still holds if one
drops the condition that $g$ is of product structure near
$\partial M$.}

$\ $

\noindent {\bf d). Some results concerning the $\eta$-invariant
associated to $g$}

$\ $

We start with Corollary 2.6 which says that the following number,
$$\left({1\over 2\pi\sqrt{-1}}\right)^{(\dim M+1)/2}\int_M
\widehat{A}\left(R^{TM}\right){\rm
Tr}\left[\exp\left(-R^E\right)\right]{\rm ch}(g,d) +
\overline{\eta}\left(\bM, g \right),\eqno(5.14)
$$ is an integer.

It is not difficult to find example such that $\overline{\eta}(\bM,
g)$ is not an integer. So it is not a trivial invariant and deserves
further study.

Our first result will show that the number in (5.14) is still an
integer if $P_{\bM}(L)$ is changed to a $Cl(1)$-spectral section in
the sense of Melrose and Piazza [MP].

Thus let $P$ be a $Cl(1)$-spectral section associated to
$P_{\partial M}(L)$. That is, $P$ defers from $P_{\partial M}(L)$
only by a finite dimensional subspace. Then $P$, as well as
$g^{-1}Pg $, is still a self-adjoint elliptic boundary condition
for $D^E$ and in view of [DF], all our previous discussion carries
over. We can thus define the corresponding $\eta$-invariant
$\overline{\eta}(\bM, g, P)$ similarly (we inserted $P$ in the
notation to emphasize its dependence).

$\ $

\noindent {\bf Proposition 5.6.} {\it The following identity holds
for any two $Cl(1)$-spectral sections $P$, $Q$ associated to
$P_{\partial M}(L)$, }
$$\overline{\eta} \left(\bM, g, P\right)\equiv
\overline{\eta} \left(\bM, g, Q\right) \ \ \ {\rm mod}\ {\bf
Z}.\eqno(5.15)$$

{\it Proof.} Let $D^{[0, 1], \psi}_{  P,g^{-1}Pg}$ denote the
elliptic self adjoint operator defined by $ D^{\psi,g}$ on $[0,
1]\times \bM$, with the boundary condition $P   $ on $\{0\}\times
\bM$ and the boundary condition ${\rm Id}-g^{-1}Pg  $ at
$\{1\}\times \bM$. Then we have
$$
\overline{\eta}\left( \bM, g, P \right) - \overline{\eta}\left(
\bM, g, Q \right)\equiv \overline{\eta} \left(D^{[0,1], \psi}_{P
 ,g^{-1}Pg}\right)- \overline{\eta} \left(D^{[0,1], \psi}_{
Q ,g^{-1}Qg} \right) \ \ \ {\rm mod}\ {\bf Z} .\eqno(5.16)
$$

In view of the definition of the operator $D^{[0,1], \psi}$ and
using the mod ${\bf Z}$ version of [KL, Theorem 7.7] repeatedly,
one deduces that
$$
\overline{\eta} \left(D^{[0,1], \psi}_{P ,g^{-1}Pg}\right)-
\overline{\eta} \left(D^{[0,1], \psi}_{ Q,g^{-1}Qg} \right)
 $$
$$
 =\overline{\eta} \left(D^{[0,1], \psi}_{ P,g^{-1}Pg}\right)-
\overline{\eta} \left(D^{[0,1], \psi}_{Q,g^{-1}Pg   } \right)+
\overline{\eta} \left(D^{[0,1], \psi}_{Q, g^{-1}Pg}\right)-
\overline{\eta} \left(D^{[0,1], \psi}_{Q ,g^{-1}Qg } \right)
$$
$$
\equiv \overline{\eta} \left(D^{[0,1]
}_{P,Q}\right)-\overline{\eta} \left((g^{-1}Dg)^{[0,1]
 }_{g^{-1}Pg, g^{-1}Qg} \right) \ \ \ {\rm mod}\ {\bf Z}.
\eqno(5.17)
$$

On the other hand, on clearly has
$$
\overline{\eta} \left(D^{[0,1] }_{P,Q}\right) \equiv
\overline{\eta} \left((g^{-1}Dg)^{[0,1] }_{g^{-1}Pg, g^{-1}Qg}
\right) \ \ \ {\rm mod}\ {\bf Z}.\eqno(5.18)$$

From (5.17) and (5.18), we obtain (5.16).  Q.E.D.

$\ $

By Proposition 5.6, when mod ${\bf Z}$, $\overline{\eta}(\bM, g, P)$
depends only on $D^E_{\partial M}$ and $g|_{\partial M}$. From now
on we denote this ${\bf R}/{\bf Z}$-valued function by
$\overline{\eta} (D^E_{\partial M},g)$.

$\ $

\noindent{\bf Remark 5.7.} In fact, for any closed spin manifold $X$
of even dimension, if the canonical Dirac operator $D^E_X$ has
vanishing index (then by [MP] there exist the associated
$Cl(1)$-spectral sections), one can define $\overline{\eta}
(D^E_{X},g)$ for $g:X\rightarrow U(N)$.

$\ $

The next  result describes the
 dependence of
$\overline{\eta} (D^E_{\partial M},g)$ on $g|_{\partial M}$.

$\ $

\noindent {\bf Theorem 5.8.} {\it If $\{g_t\}_{0\leq t\leq 1}$ is a
smooth family of maps from $M$ to $U(N)$, then
$$\overline{\eta} \left(D^E_{\partial M},g_1\right)-\overline{\eta}
\left(D^E_{\partial M},g_0\right)$$
$$\equiv -\left({1\over 2\pi\sqrt{-1}}\right)^{\dim M+1\over 2}\int_{\partial M}
\widehat{A}\left(R^{TM}\right){\rm
Tr}\left[\exp\left(-R^E\right)\right] \widetilde{\rm ch}(g_t,d,0\leq
t\leq 1)\ \  {\rm mod}\ {\bf Z}.\eqno(5.18)$$ In particular, if
$g_0={\rm Id}$, that is, $g=g_1$ is homotopic to the identity map,
then}
$$\overline{\eta} \left(D^E_{\partial M},g\right) $$
$$\equiv -\left({1\over 2\pi\sqrt{-1}}\right)^{\dim M+1\over 2}\int_{\partial M}
\widehat{A}\left(R^{TM}\right){\rm
Tr}\left[\exp\left(-R^E\right)\right] \widetilde{\rm ch}(g_t,d,0\leq
t\leq 1)\ \  {\rm mod}\ {\bf Z}.\eqno(5.19)$$

{\it Proof.} By the integrality of the number in (5.14), Proposition
5.5 and the definition of $\overline{\eta} \left(D^E_{\partial
M},g_t\right)$, one finds
$$
\overline{\eta} \left(D^E_{\partial M},g_t\right) \equiv
-\left({1\over 2\pi\sqrt{-1}}\right)^{\dim M+1\over 2}\int_M
\widehat{A}\left(R^{TM}\right){\rm
Tr}\left[\exp\left(-R^E\right)\right] {\rm ch}(g_t,d)\ \ {\rm
mod}\ {\bf Z}.\eqno(5.20)
$$
By Theorem 5.4 and (5.20) one deduces that
$$
\overline{\eta} \left(D^E_{\partial M},g_1\right)-\overline{\eta}
\left(D^E_{\partial M},g_0\right)
$$
$$
\equiv- \left({1\over 2\pi\sqrt{-1}}\right)^{\dim M+1\over
2}\int_M \widehat{A}\left(R^{TM}\right){\rm
Tr}\left[\exp\left(-R^E\right)\right] d\widetilde{\rm
ch}(g_t,d,0\leq t\leq 1)\ \  {\rm mod}\ {\bf Z}
$$
$$
\equiv- \left({1\over 2\pi\sqrt{-1}}\right)^{\dim M+1\over
2}\int_{\partial M} \widehat{A}\left(R^{TM}\right){\rm
Tr}\left[\exp\left(-R^E\right)\right] \widetilde{\rm
ch}(g_t,d,0\leq t\leq 1)\ \ {\rm mod}\ {\bf Z},\eqno(5.21)
$$ which is exactly
(5.18).

(5.19) follows from (5.18) immediately. Q.E.D.

$\ $

\noindent {\bf Remark 5.9.} As we mentioned in Remark 2.5,
$\overline{\eta}(\bM, g)$ gives an intrinsic interpretation of the
Wess-Zumino term in the WZW theory. When $\bM=S^2$, the Bott
periodicity tells us that every $K^1$ element $g$ on $S^2$ can be
deformed to the identity (adding a trivial bundle if necessary).
Hence, (5.19) gives another intrinsic form of the Wess-Zumino
term, which is purely local on $S^2$.

$\ $

Now let $\widetilde{g}^{TM}$ (resp. $(\widetilde{g}^E,\widetilde{\nabla}^E)$)
be  another Riemannian metric (resp. another couple of Hermitian metric and
connection) on $TM$ (resp. $E$). Let $\widetilde{R}^{TM}$
(resp. $\widetilde{R}^E$) be the curvature of $\widetilde{\nabla}^{TM}$
(resp. $\widetilde{\nabla}^E$), the Levi-Civita
connection of $\widetilde{g}^{TM}$.

Let $\widetilde{D}^E$ be the corresponding (twisted) Dirac
operator.

Let $\omega$ be the Chern-Simons form which transgresses the  $\widehat{A}\wedge {\rm ch}$ forms:
$$
d\omega= \left({1\over 2\pi\sqrt{-1}}\right)^{\dim M+1\over 2}
\left(\widehat{A}\left(\widetilde{R}^{TM}\right)\left[\exp\left(-
\widetilde{R}^E\right)\right] -
\widehat{A}\left({R}^{TM}\right)\left[\exp\left(-
{R}^E\right)\right]\right).\eqno(5.22)
$$

One then has the following formula describing the variation of
$\overline{\eta} \left(D^E_{\partial M},g\right)$, when
$g^{TM}|_{\partial M}$, $g^E|_{\partial M}$ and $\nabla^E|_{\partial M}$
change.

$\ $

\noindent {\bf Theorem 5.10.} {\it The following identity holds,}
$$
\overline{\eta} \left(\widetilde{D}^E_{\partial M},g\right)-
\overline{\eta} \left(D^E_{\partial M},g\right) \equiv-
\int_{\partial M} \omega {\rm ch}(g,d)\ \  {\rm mod}\ {\bf
Z}.\eqno(5.23)
$$

{\it Proof.} The proof of (5.23) follows directly from (5.22) and
the integrality of the numbers of form (5.14). Q.E.D.

$\ $

As the last result of this subsection, we prove an additivity
formula for $\overline{\eta} (D^E_{\partial M},g)$.

$\ $

\noindent{\bf Theorem 5.11.} {\it Given $f$, $g: M\rightarrow U(N)$,
the following identity holds in ${\bf R/Z}$,}
$$
\overline{\eta} \left(D^E_{\partial M},fg\right)=
\overline{\eta} \left(D^E_{\partial M},f\right) +
\overline{\eta} \left(D^E_{\partial M},g\right) .\eqno(5.24)
$$

{\it Proof.} Let $P$ be a $cl(1)$-spectral section for
$D^E_{\partial M}$ in the sense of [MP]. By [KL, Theorem 7.7], one
deduces that in ${\bf R/Z}$,
$$
\overline{\eta} \left(D^E_{\partial M},fg\right)= \overline{\eta}
\left(D^E_{\partial M},f\right) + \overline{\eta} \left((f^{-1}D^E
f)_{\partial M},fg\right) .\eqno(5.25)
$$

On the other hand, by proceeding as in (5.18), one deduces that the
following formula holds in ${\bf R/Z}$,
$$
\overline{\eta} \left((f^{-1}D^Ef)_{\partial M},fg\right)  =
 \overline{\eta} \left(D^E_{\partial M},g\right)
.\eqno(5.26)
$$
From (5.25) and (5.26), (5.24) follows.  Q.E.D.

$\ $

\noindent {\bf Remark 5.12.} Formulas (5.18), (5.19), (5.23) and
(5.24)   still hold if $\partial M$ is replaced by a closed even
dimensional spin manifold $X$ on which the Dirac operator $D^E_X$
has vanishing index, and $g:X\rightarrow U(N)$ is defined only on
$X$.

$\ $

\noindent {\bf Remark 5.13.} It might be interesting to note the
duality that $\overline{\eta} (D^E_{\partial M},g)$
 is a spectral invariant associated to a $K^1$-representative
on an {\it even} dimensional manifold, while the usual Atiyah-Patodi-Singer $\eta$-invariant
([APS1]) is a spectral invariant associated to a $K^0$-representative on an {\it odd}
dimensional manifold.

$\ $

\section*{Appendix: Toeplitz index and the Atiyah-Patodi-Singer index theorem}

In this appendix we outline a new proof of (1.1), which computes the
index of Toeplitz operators on closed manifolds. We use the notation
in Section 2, but we assume instead that the odd dimensional
manifold $M$ has no boundary.

We form the cylinder $[0,1]\times M$ and pull back everything to it
from $M$. We also identify $S(TM)$ with $S_+(T([0,1]\times
M))|_{\{i\}\times M}$, $i=0,\ 1$. Let $\widetilde{D}^E$ now be the
twisted Dirac operator on $[0,1]\times M$ acting on
$\Gamma(S_+(T([0,1]\times M))\otimes E\otimes {\bf C}^N)$. Then
$$P^E:L^2\left(S(TM)\otimes E\otimes {\bf C}^N\right)\rightarrow L^2_{\geq 0}\left(S(TM)\otimes E\otimes {\bf C}^N\right)$$
is the Atiyah-Patodi-Singer boundary condition for $\widetilde{D}^E$ at $\{ 0\}\times M$.
We equip the generalized Atiyah-Patodi-Singer boundary condition ${\rm Id}-gP^Eg^{-1}$ at $\{1\}\times M$.

Let $\widetilde{D}^E_{P^E,gP^Eg^{-1}}$ denote the elliptic operator with the Atiyah-Patodi-Singer boundary condition at
$\{ 0\}\times M$ and with the boundary condition ${\rm Id}-gP^Eg^{-1}$ at $\{1\}\times M$.

By using the standard variation formula for the index of elliptic boundary problems of Dirac type operators
(cf. [BW]), one deduces directly that

$\ $

\noindent{\bf Theorem A.1} {\it The following identity holds,}
$${\rm ind}\, T^E_g={\rm ind}\, \widetilde{D}^E_{P^E,gP^Eg^{-1}}.\eqno(A.1)$$

$\ $

Now let $\phi: [0,1]\rightarrow [0,1]$ be an increasing function such that $\phi(u)=0$ for $0\leq u\leq {1\over 4}$ and
$\phi(u)=1$ for ${3\over 4}\leq u\leq 1$. Morever, let $\widehat{D}^E$ be the Dirac type operator on
$[0,1]\times M$ such that for any $0\leq u\leq 1$,
$$\widehat{D}^E(u)=(1-\phi(u))\widetilde{D}^E +\phi(u)g\widetilde{D}^Eg^{-1}.\eqno(A.2)$$

Let $\widehat{D}^E_{P^E,gP^Eg^{-1}}$ denote the elliptic boundary value problem for $\widehat{D}^E$ with
the boundary condition $P^E$ at
$\{ 0\}\times M$ and with the boundary condition ${\rm Id}-gP^Eg^{-1}$ at $\{1\}\times M$. Then by the homotopy
invariance of the index of Fredholm operators, one has directly that
$${\rm ind}\, \widetilde{D}^E_{P^E,gP^Eg^{-1}}={\rm ind}\,
\widehat{D}^E_{P^E,gP^Eg^{-1}}.\eqno(A.3)$$

Now one can apply the Atiyah-Patodi-Singer index theorem [APS1], combined with the local index computation
 involving the Dirac type operator $\widehat{D}^E$, to get that
$${\rm ind}\,
\widehat{D}^E_{P^E,gP^Eg^{-1}}=-\left\langle \widehat{A}(TM){\rm ch}(E){\rm ch}(g),[M]\right\rangle
-\overline{\eta}\left(D^E\right)+\overline{\eta}\left(gD^Eg^{-1}\right)$$
$$=-\left\langle \widehat{A}(TM){\rm ch}(E){\rm ch}(g),[M]\right\rangle.\eqno(A.4)$$

{}From (A.1), (A.3) and (A.4), one gets (1.1).  Q.E.D.

$\ $

\noindent {\bf Remark A.2.} In view of the above proof of (1.1), one
may think of Theorem 2.3 as an index theorem on manifolds with
corners.

\noindent  Xianzhe Dai, Department of Mathematics, University of  California,
Santa Barbara, California 93106, USA

\noindent {\it E-mail}: dai@math.ucsb.edu

$\ $

\noindent Weiping Zhang, Chern Institute of Mathematics \& LPMC,
Nankai University, Tianjin 300071, P. R. China.

\noindent {\it E-mail}: weiping@nankai.edu.cn

\end{document}